\documentclass[11pt]{article}
\usepackage{amsmath}
\usepackage{amssymb}

\renewcommand{\Bbb}{\mathbb}

\newcommand{\C}{{\Bbb  C}}
\newcommand{\D}{{\Bbb D}}

\newcommand{\R}{{\Bbb  R}}

\newcommand{\N}{{\Bbb  N}}

\newcommand{\T}{{\Bbb  T}}

\renewcommand{\O}{{\cal O}}

\newcommand{\PSH}{{\operatorname{PSH}}}

\newcommand{\pr}{{\operatorname{pr}}}

\numberwithin{equation}{section}

\newtheorem{theorem+}           {Theorem}      [section]
\newtheorem{definition+}  [theorem+]  {Definition}
\newtheorem{lemma+}  [theorem+]  {Lemma}
\newtheorem{corollary+}  [theorem+]  {Corollary}
\newtheorem{proposition+}  [theorem+]  {Proposition}
\newtheorem{example+}  [theorem+]  {Example}
\newtheorem{question+}  [theorem+]  {Question}

\newenvironment{theorem}{\begin{theorem+}\sl}{\end{theorem+}\rm}

\newenvironment{lemma}{\begin{lemma+}\sl}{\end{lemma+}\rm}
\newenvironment{corollary}{\begin{corollary+}\sl}{\end{corollary+}\rm}
\newenvironment{proposition}{\begin{proposition+}\sl}{\end{proposition+}\rm}
\newenvironment{example}{\begin{example+}\rm}{\end{example+}\rm}
\newenvironment{question}{\begin{question+}\rm}{\end{question+}\rm}
\newenvironment{proof}{\medbreak\noindent{\it Proof:}\rm}{\hfill$\square$\rm}
\newenvironment{prooftx}[1]{\medbreak\noindent{\it #1:}\rm}{\hfill$\square$\rm}

\title{\Large \bf THE RELATIVE EXTREMAL FUNCTION FOR BOREL
SETS IN COMPLEX MANIFOLDS}
\author{\large\bf  Armen Edigarian\footnote{This work is a part of 
the Research Grant No.~1 PO3A 005 28, which is
supported by public means in the programme for promoting science in
Poland in the years 2005-2008. The first author is a fellow of
Krzy\.zanowski Fund at the Jagiellonian University.  This work was
also supported by the Research Fund at the University of Iceland.}
\  and  Ragnar Sigurdsson }

\date{July 12, 2006}

\begin{document}
\maketitle

\begin{abstract}  \noindent
We study a disc formula for the relative extremal function for
Borel sets in complex manifolds.


\end{abstract}

\section{Introduction}\label{sec:intro}

Let $X$ be a complex manifold, $\PSH(X)$ be the class of all
plurisubharmonic functions on $X$, and $\PSH^-(X)$ be the subclass
of all non-positive functions.  For any subset $A$ of $X$
we define
$$
\omega(\cdot,A,X)=\sup\{u\in \PSH^-(X)\,;\, u|A\leq -1\}
=\sup\{u\in \PSH(X)\,;\, u\leq -\chi_A\},
$$
where $\chi_A$ denotes the characteristic function of the set $A$.
The least upper semicontinuous majorant $\omega^*(\cdot,A,X)$ of
$\omega(\cdot,A,X)$ is plurisubharmonic and it is called the
{\it relative extremal function for $A$ in $X$}.

Observe that if $A$ is a Borel set,
$u\in \PSH(X)$,  $u\leq -\chi_A$, $x=f(0)$, where
$f\in \O(\D,X)\cap C(\overline \D,X)$,
i.e., $f$ is an {\it analytic disc}  which extends to a continuous map
from the closure $\overline \D$ of the unit disc $\D$
to $X$, then the subaverage property of $u$ implies 
$$
u(x)\leq \dfrac 1{2\pi}\int_0^{2\pi}u(f(e^{i\theta}))\, d\theta
\leq \dfrac 1{2\pi}\int_0^{2\pi}(-\chi_A)(f(e^{i\theta}))\, d\theta
=-\sigma(f^{-1}(A)\cap \T)=-\sigma_f(A),
$$
where $\sigma$ denotes the normalized arc length measure on the unit
circle $\T$ and $\sigma_f$  the image measure (push-forward)
of $\sigma$ under the map $f$.   By taking  supremum over all
plurisubharmonic  $u\leq -\chi_A$ and infimum
over all $f\in \O(\D,X)\cap C(\overline \D,X)$, we 
get 
\begin{align*}
\omega(x,A,X)\leq \Omega(x,A,X)&=\inf \{-\sigma_f(A)\,;\, f\in \O(\D,X)\cap C(\overline
\D,X), f(0)=x\}\\
&=-\sup\{\sigma_f(A)\,;\, f\in \O(\D,X)\cap C(\overline
\D,X), f(0)=x\}.
\end{align*}

In this paper we are mainly concerned with a possible converse of this
inequality.  If $A$ is an open subset
of $X$, then  $\omega(\cdot,A,X)=\Omega(\cdot,A,X)$.
This is a special case of Poletsky's theorem which states that
if $\varphi$ is an upper semicontinuous function on $X$, then
for every $x$ in $X$
$$
\sup\{u(x)\,;\, u\in \PSH(X), u\leq \varphi\}
=\inf\{\int_\T \varphi\circ f\, d\sigma \,;\, f\in \O(\overline\D,X),
f(0)=x\}.
$$
See \cite{LarSig1}, \cite{LarSig3}, \cite{P0}, and \cite{R}.  Here
$\O(\overline\D,X)\subset \O(\D,X)\cap C(\overline \D,X)$ denotes the 
set of all {\it closed} analytic discs in $X$, i.e., analytic discs
which extend to holomorphic maps in some neighbourhood of 
$\overline \D$.
With $-\chi_A$ in the role of $\varphi$ we get 
$\omega(\cdot,A,X)=\Omega(\cdot,A,X)$ for every open
set $A$.

We say that the subset $A$ of $X$ is {\it pluriregular} at the point $x\in A$ if
$\omega^*(x,A,X)=-1$, we say that $A$ is {\it locally pluriregular
} at  the point $x$ in $A$ if $\omega^*(x,A\cap U,U)=-1$ for every
neighbourhood $U$ of $x$, and finally we say that $A$ is ({\it
locally}) {\it pluriregular} if $A$ is  (locally) 
pluriregular at each of its points.  Note that if $A$ is locally
pluriregular, then $A$ is pluriregular and that if $A$ is
pluriregular, then $\omega(\cdot,A,X)=\omega^*(\cdot,A,X)$.
Our main result of  Section 2
is that if $A$ is a locally pluriregular subset of $X$, then
$\Omega(\cdot,\overline A,X)\leq \omega(\cdot,A,X)$.  (See
Th.~\ref{th:1}.) This is a generalization of Th.~7.2 in
Poletsky \cite{P0}.

Let $E$ be a subset of a complex manifold $X$. We say that $E$ is
{\it pluripolar} or {\it locally pluripolar} if for any $a\in E$
there exists a neighbourhood $U$ of $a$ in $X$ and $u\in \PSH(U)$,
$u\not\equiv-\infty$, such that $E\cap U\subset\{x\in U \,;\,
u(x)=-\infty\}$. We say that $E\subset X$ is  {\it globally
pluripolar} if there exists $u\in\PSH(X)$, $u\not\equiv-\infty$,
with $E\subset \{x\in X; u(x)=-\infty\}$. Note that any globally
pluripolar set is locally pluripolar. Josefson \cite{J} proved
that in $\C^n$ every pluripolar set is globally pluripolar.
We say that a complex manifold $X$ is a {\it Josefson manifold} if
any locally  pluripolar set is globally pluripolar. Bedford
\cite{Be} has generalized Josefson's theorem to a certain class of
complex spaces including Stein manifolds.  He also showed that
examples, originally given by Grauert \cite{Gr}, of complex manifolds
which possess no non-constant holomorphic functions  are Josefson
manifolds.

In Section 3 we prove that if $X$ is a relatively compact domain
in a Josefson manifold and $A$ is a Borel subset of $X$, then
$\Omega(\cdot,A,X)\leq \omega^*(\cdot,A,X)$. (See Th.~3.1.)  
The main result of Section 3 is that a
Josefson manifold $X$ has the property that every bounded
plurisubharmonic function on $X$ is constant if and only if
for every $p\in X$, every non-pluripolar Borel subset $A$ of $X$, and every $\varepsilon>0$
there exists $f\in\O(\D,X)\cap C(\overline\D, X)$ such that $f(0)=p$ and
\begin{equation*}
\sigma_f(A)=\sigma(f^{-1}(A)\cap\T)>1-\varepsilon.
\end{equation*}
(See Th.~\ref{th:2}.) As a consequence we get a
characterization of pluripolar sets in terms of analytic discs.

In Section 4 we look at Borel subsets $A$ of the boundary
${\partial}D$ of a relatively compact domain in a complex manifold
$X$.  We define the relative extremal function for an open
subset $U$ of the boundary as $\omega(\cdot,U,D)=u_{-\chi_U,D}$, where
$u_{f,D}$ is the Perron-Bremermann envelope of the boundary function 
$f$, and for any subset $A$ of ${\partial}D$ we define
$\omega(\cdot,A,D)$ as the supremum over all $\omega(\cdot,U,D)$ for
$U$ open containing $A$.  We call the domain {\it weakly regular}
if the upper semicontinuous extension $u^*_{-\chi_U,D}$ of
$u_{-\chi_U,D}$ to the closure $\overline D$ is less than or equal to
$-\chi_{U}$ on ${\partial}D$. For any Borel subset $A$ of
${\partial}D$ we define $\Omega(x,A,D)$ as the infimum over 
$-\sigma_f(A)$ for $f\in\O(\D,D)\cap C(\overline \D,\overline D)$
with $f(0)=x$.  We prove (see Th.~4.3) that for a weakly regular
domain $D$ and  every open subset $U$ of ${\partial}D$
we have $\omega(\cdot,U,D)=\Omega(\cdot,U,D)$  
and (see Th.~4.10)
$\omega(\cdot,\overline A,D)\leq \Omega(\cdot,\overline A,D)\leq
\omega^*(\cdot,A,D)$, if $A$ is a Borel subset of the form
$A=A_1\cup E$, where $A_1$ is locally pluriregular with respect to $D$
and $E$ is such that there exists $u\in \PSH^-(D)$,
$u\not\equiv-\infty$, and $u^*|_E\equiv -\infty$.  It remains an
open question if the last result holds for every Borel set $A$.

\section{Construction of analytic discs}

We have already seen that for every manifold $X$ and every Borel
subset $A$ of $X$ we have
$
\omega(\cdot,A,X)\leq \Omega(\cdot,A,X)
$
and that Poletsky's theorem implies that equality holds if $A$ is
open.

\begin{theorem}\label{th:1}
Let $X$ be a complex manifold and $A$ be any
locally pluriregular subset of $X$.  Then $\Omega(\cdot,\overline A,X)\leq
\omega(\cdot,A,X)$, in particular, $\Omega(\cdot,A,X)=\omega(\cdot,A,X)$ if $A$ is also closed. 
\end{theorem}

The main argument of the proof consists of an approximation of analytic
discs and it appears a few times in this paper.  We therefore  state
it as a separate result.    A similar result
for  domains in $\C^n$ and,  more generally, for domains in Banach spaces
is proved by   Poletsky in  \cite{P1}.  Our proof uses the existence
of Stein neighbourhoods of certain sets which was  proved by
Rosay  \cite{R}.  For a simplification of his arguments and further
development  see \cite{E2} and \cite{LarSig3}.

If $X$ is a complex manifold and $d:X\times X\to [0,+\infty)$ is
a continuous function vanishing on the diagonal,
i.e., $d(x,x)=0$ for all $x\in X$, then
for any subset $A$ of $X$ we define the {\it diameter of $A$ with
respect to $d$ } as
$\sup\{d(x,y)\,;\, x,y\in A\}$.
In the proof of Theorem \ref{th:1} we will take $d$ as a complete
hermitian metric defining the topology of $X$.

\begin{theorem}\label{th:family}
Let  $X$ be a complex manifold, $d:X\times X\to [0,+\infty)$ be a
continuous function vanishing on the diagonal, $\delta>0$, and
 $\{B_j\}$ be a countable family of open subsets in $X$  of
diameter less than $\delta$ with respect to $d$.   Assume that
$U$ and $V$ are open subsets of $X$ and
$$
V\subset \cup_j \{x\in B_j \,;\, \omega(x,U\cap B_j,B_j)<-a\},
$$
where $a\in (0,1)$.
Let $h\in\O(\overline \D,X)$ and assume that  $\Delta\subset h^{-1}(V)\cap \T$
is a non-empty open set.
Then for every $\varepsilon\in (0,1)$ there exist $g\in \O(\overline \D,X)$
and an open set $\tilde \Delta\subseteq \Delta$  such that
\begin{description}
\item{(1)} $g(0)=h(0)$,
\item{(2)} $d(g,h)=\sup_{z\in \overline\D}d(g(z),h(z))\leq \delta+\varepsilon$,
\item{(3)} $\sigma(\tilde\Delta)\geq (1-\varepsilon)a\sigma(\Delta)$,
and
\item{(4)} $g(\tilde\Delta)\subset U$.
\end{description}
\end{theorem}

\begin{proof}
For $r>0$ we let  $D_r$ be the open disc in $\C$ with radius $r$ 
and centre at the origin and we assume that $h\in\O(D_s, X)$ for some $s>1$.
Fix $\Delta_0\subset\Delta$ a union of closed arcs such that 
$\sigma(\Delta_0)>(1-\varepsilon)\sigma(\Delta)$.

Take $w_0\in\Delta$. Then $x_0=h(w_0)\in V$,  so there exists a $j_0$
such that  $x_0\in B_{j_0}$ and
$\omega(x_0,U\cap B_{j_0},B_{j_0})<-a$.  Since $U$ is open,
Poletsky's theorem implies that
there exists  $f_0\in\O(\overline\D,B_{j_0})$
such that $f_0(0)=x_0$ and $\sigma_{f_0}(U)>a$.
Let $I_0\subset f_0^{-1}(U)\cap\T$ be a union of
finite number of closed arcs such that $\sigma(I_0)>a$.
By Lemma 2.3 in \cite{LarSig1}, there
exists an open neighbourhood $V_0$ of $x_0=f_0(0)$ in $X$,
$r>1$, and $f\in\O(D_r\times V_0,B_{j_0})$ such that
$f(z,x_0)=f_0(z)$ for all $z\in D_r$ and  $f(0,x)=x$ for all $x\in V_0$.
By choosing $r>1$ sufficiently small and shrinking the neighbourhood $V_0$ of
$x_0$,  we may assume that  $f(z,x)\in U$ for all $z\in I_0$ and $x\in V_0$.
We set $F_0(z,w)=f(z,h(w))$ and note that $F_0$ is defined
on $D_r\times h^{-1}(V_0)$ and that $h^{-1}(V_0)$ is a neighbourhood of $w_0$.

We apply  a compactness argument on $\Delta_0$ and conclude that we may find:
\begin{description}
\item{$\bullet$} Open discs $U_1,\dots,U_m$  centred on $\T$ with
mutually disjoint closures such that $U_j\cap\T\subset\Delta$ and
$\sigma((U_1\cup\dots \cup
U_m)\cap\T)>(1-\varepsilon)\sigma(\Delta)$. %
\item{$\bullet$} $r_j>1$, $j=1,\dots,m$, and  
holomorphic maps $F_j:D_{r_j}\times
U_j\to B_{k(j)}$ with  $F_j(0,w)=h(w)$ for all
$w\in U_j$. 
\item{$\bullet$} Finite unions $I_j$, $j=1,\dots,m$, of closed arcs 
  on $\T$ with $\sigma(I_j)>a$ and
$F_j(z,w)\in U$ for all $z\in I_j$ and $w\in U_j$.
\end{description}

Take closed arcs $J_1,\dots,J_m$ in $\T$ such that $J_j\subset
U_j\cap \T$ and 
$\sigma(J_1\cup\dots\cup J_m)>(1-\varepsilon)\sigma(\Delta)$.
Let
\begin{equation*}
K_0=\{(w,0,0,0,h(w))\,;\, w\in\overline\D\} \subset \C^4\times X
\end{equation*}
and
\begin{equation*}
K_j=\{(w,z,0,0,F_j(z,w))\,;\, w\in J_j, z\in\overline\D\} \subset
\C^4\times X,
\qquad j=1,\dots,m.
\end{equation*}

By the proof of Th.~1.2 in \cite{LarSig3}, there exists a
Stein neighbourhood  $Z$ of $K_0\cup K_1\cup\dots \cup K_m$ in
$\C^4\times X$. Let $\tau:Z\to\C^N$ be an embedding,
$\varkappa:W\to \tau(Z)$  be a holomorphic retraction from a Stein
neighbourhood $W$ of $\tau(Z)$ in $\C^N$, and  $\varphi=\pr\circ
\tau^{-1}\circ\varkappa:W\to X$ be the holomorphic submersion,
where $\pr:\C^4\times X\to X$ is the projection.

We let $\rho:\T\to[0,1]$ be a $C^\infty$ function such that
$\rho=0$ on $\T\setminus (\cup_j J_j)$ and
$\rho=1$ on a subset of $\cup_j J_j$ such that
$\sigma(\{w\in\T \,;\,  \rho(w)=1\})>(1-\varepsilon)\sigma(\Delta)$.
We define a $C^\infty$ map $F:D_s\times\T\to X$ by
$$
F(z,w)=
\begin{cases}
 F_j(\rho(w)z,w),&\quad w\in J_j, \ j=1,\dots,m,\\
 h(w),&\quad w\not\in\cup_j J_j.
\end{cases}
$$
Since $Z$ is a neighbourhood of $K_0\cup\cdots\cup K_m$
we can  replace $s>1$ by a smaller number and
can define $\widetilde F:D_s\times \T\to W$ by
 $\widetilde F(z,w)=\tau(w,z,0,0,F(z,w))$ and
$\widetilde h:D_s\to W$ by $\widetilde
h(w)=\tau(w,0,0,0,h(w))$.  We note that
$\widetilde F(0,w)=\widetilde h(w)$ for all $w\in\T$.

In exactly the same way as in the proof of Lemma 2.6 in
\cite{LarSig1} we construct a sequence $\widetilde
F_j\in\O(D_s\times A_j,W)$, $j\ge j_0$, where $A_j$ is an open
annulus containing $\T$, such that
\begin{description}
\item{$\bullet$} $\widetilde F_j\to\widetilde F$ uniformly on
$D_s\times\T$ as $j\to\infty$, %
\item{$\bullet$} there is an integer $k_j\ge j$ such that 
for all $k\geq k_j$ the map
$\widetilde F_j(zw^k,w)$ can be extended to a map $\widetilde
G_j\in\O( D_{s_j}\times D_{s_j},W)$, where $s_j\in(1,s)$, and
\item{$\bullet$} $\widetilde G_j(0,w)=\widetilde h(w)$ for all
$w\in D_{s_j}$.
\end{description}

We need to estimate
$\sup_{w\in\T}\sup_{z\in\overline\D}d(\varphi(\widetilde F_j(z,w)),h(w))$.
Since $\widetilde F_j\to\widetilde F$ uniformly on $\overline \D\times\T$ as
$j\to\infty$ and $\varphi(\tilde F)=F$ we have
$$
\sup_{z\in\overline\D, w\in\T}d(\varphi(\widetilde F_j(z,w)),h(w))\to
\sup_{z\in\overline\D,w\in\T}d( F(z,w),h(w)), \qquad j\to\infty.
$$
We have $\sup_{z\in\overline\D}d(F(z,w),h(w))=0$
for all  $w\in\T\setminus(\cup_j J_j)$ and since $F_j$ takes values 
in $B_{k(j)}$, we have
$\sup_{z\in\overline\D}d(F(z,w),h(w))\le \delta$
for all  $w\in J_j$.  Hence
$$
\limsup_{j\to\infty}\Big[\sup_{z\in\overline\D,w\in\T}d(\varphi(\widetilde
F_j(z,w)),h(w))\Big]\le \delta.
$$

Take $j\ge j_0$ so that $\sup_{z\in\overline\D,w\in\T}
d(\varphi(\widetilde F_j(z,w)),h(w))\le \delta+\varepsilon/2$.
There exists $t\in(0,1)$ such that
$\sup_{z\in\overline\D, w\in[t,1]}d(\varphi(\widetilde
F_j(z,w)),h(w))< \delta+\varepsilon$.
Note that
$$
\sup_{z\in\overline\D,|w|\le t}d(\varphi(\widetilde F_j(zw^k,w)),h(w))\to0
$$
as $k\to\infty$,  so for sufficiently large $k$ we have
$$
\sup_{z,w\in\overline\D}d(\varphi(\widetilde F_j(zw^k,w)),h(w))< \delta+\varepsilon.
$$
We set $G(z,w)=\varphi(\widetilde F_j(zw^k,w))$. Then
$G\in\O(\overline\D^2,X)$ and $G(0,w)=h(w)$ for all $w\in\D$.

Put $C=\cup_j (I_j\times \widetilde J_j)$, where $\widetilde J_j=\{w\in J_j\,;\, \rho(w)=1\}$.
If $\sigma_2=\sigma\times \sigma$ is the product measure on 
the torus $\T^2$, then
$$\sigma_2(C)=\sum_j\sigma_2(I_j\times \widetilde J_j)=
\sum_j\sigma(I_j)\sigma(\widetilde J_j)>
a(1-\varepsilon)\sigma(\Delta).
$$
The map 
$\T^2\ni(z,w)\to (zw^k,w)\in\T^2$ is an automorphism with the absolute
value of the Jacobian equal to $1$. Therefore the measure
of the set $\widetilde C=\{(z,w)\in\T^2\,;\,  (zw^k,w)\in C\}$ is equal to $\sigma_2(C)$. By Fubini's
theorem there is a $\theta\in[0,2\pi)$ such that $\sigma(C')\ge\sigma_2(\widetilde C)$, where
$C'=\{w\in\T\,;\, (e^{i\theta}w,w)\in\widetilde C\}$.

Now we finally define  $g(w)=G(e^{i\theta}w,w)$ for $w\in \overline
\D$ and
$\widetilde\Delta=g^{-1}(U)\cap\Delta$. 
Then (1) and (4) are obvious and (2) holds because
$$
d(g,h)\le\sup_{z,w\in\overline\D}d(G(z,w),h(w))<\delta+\varepsilon.
$$
For proving (3) we take  $w\in C'$ and observe that 
 $(e^{i\theta}w,w)\in \widetilde C$ and therefore $(e^{i\theta}w\cdot w^k,w)\in C$.
This implies that $e^{i\theta}w^{k+1}\in I_j$,
$w\in\widetilde J_j$ for some $j$, and consequently $g(w)\in U$. 
Hence
$C'\subset\widetilde\Delta$ and 
$$\sigma(\widetilde\Delta)\geq \sigma(C')\geq \sigma_2(C)>(1-\varepsilon)a\sigma(\Delta).
$$
\end{proof}

\begin{prooftx}{Proof of Th.~\ref{th:1}}
Let $x_0\in X$.  It is sufficient to prove that if $a\in (0,1)$ and
$\omega(x_0,A,X)<-a$, then $\Omega(x_0,\overline A,X)\leq -a$.  This
inequality will in turn follow if we prove that for every
$\varepsilon\in (0,1)$ there exists $h\in \O(\D,X)\cap
C(\overline\D,X)$ such that $h(0)=x_0$ and 
$\sigma_h(\overline A)>(1-\varepsilon)a$.

We take $\varepsilon_m\searrow 0$ such that
$\prod_m (1-\varepsilon_m)\geq
\sqrt{1-\varepsilon}$.  For every $m$ we find a covering
$\{B^m_j\}$ of $X$ by countably many balls of diameter less than $\varepsilon_m$ and set
$$
U_m=\cup_j \{x\in B_j^m\,;\, \omega^*(x,A\cap
B_j^m,B_j^m)<-1+\varepsilon_m\}.
$$
Since $A$ is locally pluriregular, $U_m$ is a neighbourhood of $A$
and the inequality 
$$\omega(\cdot,A\cap B_j^m,B_j^m)\geq
\omega(\cdot,U_{m+1}\cap B_j^m,B_j^m)$$ 
implies
$$
U_m\subseteq
\cup_j \{x\in B_j^m\,;\, \omega(x,U_{m+1}\cap
B_j^m,B_j^m)<-1+\varepsilon_m\}.
$$
Since $A$ is locally pluriregular and $U_1$ is an open neighbourhood
of $A$ we have
$$
-a\geq \omega^*(x_0,A,X)\geq \omega(x_0,U_1,X)=\Omega(x_0,U_1,X)
$$
and there exists $h_1\in \O(\overline\D,X)$ such that $h_1(0)=x_0$
and $\sigma_{h_1}(U_1)>a$.  We set $\Delta_1=h_1^{-1}(U_1)\cap \T$
and observe that by the definition of the measure $\sigma_{h_1}$
we have $\sigma(\Delta_1)>a$.  We apply Th.~\ref{th:family} and get
inductively a sequence $h_m$ in $\O(\overline \D,X)$ and a
decreasing sequence $\Delta_m$ of open subsets of $\T$ such that
 $h_m(0)=x_0$,  $h_m(\Delta_m)\subset U_m$,
$\sigma(\Delta_{m+1})>(1-\varepsilon_m)^2\sigma(\Delta_m)$, and
$d(h_{m+1},h_m)<2\varepsilon_m$.

The last condition  implies that $h_m$ converges uniformly on $\overline\D$
to some $h\in\O(\D,X)\cap C(\overline\D,X)$.
We set $\Delta=\cap_m\Delta_m$.  Since
$h_m(\Delta_m)\subset U_m$ and the points of $U_m$ are at a distance
less than or equal to $\varepsilon_m$ from $A$, we have
$h(\Delta)\subset \overline A$ and since
$\sigma(\Delta_{m+1})>(1-\varepsilon_m)^2\sigma(\Delta_m)$ we get
$$
\sigma_h(\overline A)\geq \sigma(\Delta)>\prod_m
(1-\varepsilon_m)^2\sigma(\Delta_1)>(1-\varepsilon)a.
$$
\end{prooftx}

\section{Characterization of pluripolar sets}\label{sec:Josefson}

Let $X$ be a complex manifold. We say that $X$ is a Josefson manifold
if any  locally pluripolar subset of $X$ is globally pluripolar.
Note that any domain in a Josefson manifold is a Josefson manifold.
In particular, any domain in $\C^n$  is a Josefson manifold.
As a direct consequence of Th.~\ref{th:1} we get
(cf. Cor.~7.2 in \cite{P0})

\begin{theorem}\label{cor:1}
Let $X$ be a relatively compact domain in a Josefson manifold
and $A$ be a Borel subset of $X$.  Then
$\Omega(\cdot,A,X)\leq \omega^*(\cdot,A,X)$.
 \end{theorem}

Before we prove the theorem  we prove the following auxiliary result.

\begin{lemma}{\rm [See Th.~8.3 in \cite{Be-Ta}
 or  Th.~7.3 in \cite{P0})]}\label{lemma:7.3}. \  Let  $\mu$ be a
Borel probability  measure which is zero on every pluripolar set.
Then the set function $c=c_\mu$
defined by
$$
c(A)=c_\mu(A)=-\int_X \omega^*(\cdot,A,X)\, d\mu,  \qquad A\subset X,
$$
is a Choquet capacity, i.e.,
\begin{description}
\item{{\rm (1)}} $c(A_1)\le c(A_2)$ if $A_1\subset A_2$; %
\item{{\rm (2)}} $c(K)=\lim_{j\to\infty} c(K_j)$, where $K_1\supset
K_2\supset\dots\supset K$ are compact sets and $K=\cap_j K_j$; %
\item{{\rm (3)}} $c(A)=\lim_{j\to\infty} c(A_j)$, where $A_1\subset
A_2\subset\dots\subset A$ are arbitrary sets and $A=\cup A_j$.
\end{description}
\end{lemma}

\begin{proof} Since $\omega^\ast(\cdot,A,X)=\omega(\cdot,A,X)$ on $X\setminus
P$ for some pluripolar set $P$, we have
$$
c(A)=c_\mu(A)=-\int_X \omega(\cdot,A,X)\, d\mu,  \qquad A\subset X.
$$
Since $-\omega(\cdot,A_1,X)\leq -\omega(\cdot,A_2,X)$ if $A_1\subset
A_2$, (1) holds.  For proving  (2) we first observe that  (1) implies 
$c(K)\le\lim_{j\to\infty} c(K_j)$. If $\{V_k\}_{k\in \N}$ is a
decreasing basis of neighbourhoods of $K$, then 
$\omega(x,V_k,X)$ increases to $\omega(x,K,X)$ and the monotone 
convergence  theorem implies that  $\lim_{k\to\infty} c(V_k)=c(K)$. 
For every $k\ge1$ there exists a $j_k$ such that $K_{j}\subset V_k$ 
for all $j\ge j_k$, so 
$$
\lim_{j\to\infty} c(K_j)\le\lim_{k\to\infty} c(V_k)=c(K).
$$
Note that  (3) is clear for open sets. Fix $\varepsilon>0$
and put $V_j=\{x\in X:\omega^\ast(x,A_j,X)<-1+\varepsilon\}$ and
$V=\cup_j V_j$. Then $c(V_j)\to c(V)$. Note that $V_j\supset
A_j\setminus P_j$, where $P_j$ is a pluripolar set. Hence,
$V\supset A\setminus P$, where $P=\cup_j P_j$.  The set $P$ 
is pluripolar, so $c(V)\ge c(A)$.

We have $\omega(\cdot,V_j,X)\ge
\omega^\ast(\cdot,A_j,X)/(1-\varepsilon)$ and, therefore, 
$c(V_j)\le c(A_j)/(1-\varepsilon)$. Hence,
$$
c(A)\le\frac{1}{1-\varepsilon} \lim_{j\to\infty} c(A_j).
$$
Since $\varepsilon>0$ is arbitrary, this proves (3).
\end{proof}

\begin{prooftx}{Proof of Th.~\ref{cor:1}}  If $A$ is pluripolar, then there
exists $u\in \PSH^-(X)$ such that  $u\not\equiv-\infty$
and $A\subset \{u=-\infty\}$. This implies that
$\omega^*(\cdot,A,X)=0$ and the inequality holds.

From  now on we assume that $A$ is non-pluripolar. Let us first
take $A$ compact. It is sufficient to show that $A$ can be written
as $A=A_1\cup E_1$, where $A_1$ is locally pluriregular and
$E_1\subset \{u=-\infty\}$ for some $u\in \PSH^-(X)$,
$u\not\equiv-\infty$. Indeed, then Theorem \ref{th:1} gives
$$
\omega^*(\cdot,A,X)=
\omega^*(\cdot,A_1,X)\geq
\Omega(\cdot,\overline A_1,X)\geq
\Omega(\cdot,\overline A,X).
$$
In order to prove that $A=A_1\cup E_1$, we choose a countable dense
subset $\{a_k\}$ of $A$ and set
$$
E_1=\bigcup_k\bigcup_m
\{x\in A\cap B_d(a_k,1/m)\,;\,
\omega^*(x,A\cap B_d(a_k,1/m),B_d(a_k,1/m))>-1\}
$$
where $B_d(a,r)$ denotes the ball with centre $a$ and radius $r$
with respect to a complete hermitian metric $d$ defining the topology
of $X$.  Note that $E_1$ is locally pluripolar and therefore by
assumption globally pluripolar.  Moreover, since $X$ is relatively
compact in a Josefson manifold, we can find $u\in\PSH^-(X)$ so that
$E_1\subset \{u=-\infty\}$.  Now we put $A_1=A\setminus E_1$.
Then $A_1$ is locally pluriregular, for if $x\in A_1$ and
$U$ is a neighbourhood of $x$, then there exists a ball
$B_d(x,1/m)\subset U$ and $a_k\in B_d(x,1/2m)$ such that $B_d(a_k,1/2m)\subset
B_d(x,1/m)$ and we get
$$
-1\leq \omega^*(x,A\cap U,U)\leq
\omega^*(x,A\cap B_d(a_k,1/2m),B_d(a_k,1/2m))=-1.
$$
Now we let $A$ be any Borel subset of $X$.  We fix $x_0\in X$ and
are going to show that
$$\Omega(x_0,A,X)\leq \omega^*(x_0,A,X).
$$
It suffices to show that there exists a sequence of compact sets
$K_1\subset K_2\subset \cdots \subset A$ so that
$$
\omega^*(x_0,K_j,X)\to \omega^*(x_0,A,X), \qquad j\to \infty.
$$
Let us construct a probability measure on $X$ which is zero
on every pluripolar set. Fix a covering $\{U_j\}$ of $X$ so that
$(U_j,\psi_j)$ is a holomorphic chart and $\psi_j(U_j)\subset\C^m$
is a bounded domain. (We assume that $X$ is $m$-dimensional.)
For any Borel set $A$ we put
$$
\mu(A)=\sum_{j=1}^\infty\frac 1{2^j}\cdot \frac{\lambda_m(\psi_j(A\cap
U_j))}{\lambda_m(\psi_j(U_j))},
$$
where $\lambda_m$ is the Lebesgue measure in $\C^m$. It is easy to
see that $\mu$ is a probability measure on $X$. Moreover, for
any pluripolar set $P$ we have $\mu(P)=0$.

By Lemma \ref{lemma:7.3} $c_\mu$ is a Choquet capacity. The Choquet
capacitability theorem states that
$$
c_\mu(A)=\sup\{c_\mu(K) \,;\, K\subset A \text{ is compact}\}
$$
for all Borel subsets $A$ of $X$. Hence, for a fixed Borel set $A$
there exists a sequence $K_1\subset K_2\subset\dots\subset A$ of
compact sets such that $c_\mu(K_j)\to c_\mu(A)$. It is easy to see
that $\omega^\ast(\cdot,K_j,X)\to\omega^\ast(\cdot, A,X)$.
\end{prooftx}

\medskip
The equivalence of (1) and (3) in the following theorem is well known
and it indeed holds on every manifold.    See  Edigarian \cite{E1} and
Rosay \cite{R}.    Using the theorem above we are able to refine this
result.

\begin{theorem}\label{th:2}
Let $X$ be a Josefson manifold. Then the following conditions are equivalent
\begin{description}
\item{(1)} Any bounded plurisubharmonic function on $X$ is constant.
\item{(2)} For every $p\in X$, every nonpluripolar Borel subset $A$
of $X$,  and every $\varepsilon>0$
there exists $f\in\O(\D,X)\cap C(\overline\D, X)$ such that $f(0)=p$ and
\begin{equation*}
\sigma_f(A)=\sigma(f^{-1}(A)\cap\T)>1-\varepsilon.
\end{equation*}
\item{(3)} For every $p\in X$, every nonempty open subset $U$ of $X$, and every $\varepsilon>0$
there exists $f\in\O(\D,X)\cap C(\overline\D, X)$ such that $f(0)=p$ and
\begin{equation*}
\sigma_f(U)>1-\varepsilon.
\end{equation*}
\end{description}
\end{theorem}

\begin{proof}
The proof that  (2) implies (3) is trivial.
In order to prove that (3) implies (1), we let
$u$ be a negative plurisubharmonic function on $X$. Assume that $u$ is non-constant.
Then there exist $x_1,x_2\in X$ such that $u(x_1)<u(x_2)$. Take an $a\in\big(u(x_1),u(x_2)\big)$.
Put $U=\{x\in X\,;\, u(x)<a\}$. Then $U$ is an open set and $x_1\in
U$.  By (3) we have  $\Omega(\cdot,U,X)\equiv-1$.
Since  $-1\leq \omega(\cdot,U,X)\le\Omega(\cdot,U,X)$, we have 
$\omega(\cdot,U,X)\equiv-1$. But $\frac{1}{|a|}
u(\cdot)\le\omega(\cdot,U,X)$, which implies  $u(x_2)\le a$, a contradiction.

In order to prove that (1) implies (2), we  take a sequence of
subdomains $X_1\Subset X_2\Subset\dots\Subset X$ such that
$\cup_{n=1}^\infty X_n=X$. There exists a compact set $K\subset A$
such that $K$ is nonpluripolar. Without loss of generality we may
assume that $K\subset X_1$. For any $n\ge1$ we have  by
Theorem~\ref{cor:1} that
\begin{equation*}
\Omega(\cdot,K,X_n)\le u_n=\omega^\ast(\cdot,K,X_n)\quad\text{ on }X_n.
\end{equation*}
There exists  $x_1\in K$ such that $u_1(x_1)=-1$. Note that the sequence $\{u_n\}$ is decreasing.
Put $u=\lim u_n\in\PSH(X)$. Hence $u$ is a constant, $u(x_1)=-1$, so $u\equiv-1$.

Fix a $p\in X$ and $\varepsilon>0$. Then $\Omega(p,K,X_n)\to-1$ as
$n\to\infty$, so there exists $n\in \N$ such that
$$
\Omega(p,F,X_n)\leq \Omega(p,K,X_n)<-1+\varepsilon.
$$
\end{proof}

\bigskip
Observe that if  $u\in \PSH(X)$ is such that $E\subseteq \{x\in
X\,;\, u(x)=-\infty\}=\tilde E$,  $x\in X\setminus \tilde E$, and
$f\in \O(\D,X)\cap C(\overline \D,X)$ with $f(0)=x$, then
$$
-\infty<u(x)\le\int_\T u\circ f\, d\sigma,
$$
and we conclude that  $\sigma_f(\tilde E)=\sigma(\{t\in \T\,;\,
u(f(t))=-\infty\})=0$.  Hence we have

\begin{theorem}\label{th:3}
Let $X$ be a complex  manifold and let $E$ be a globally
pluripolar subset.  Then there exists a globally pluripolar
$\tilde E\supseteq E$ such that for every $x\in X\setminus \tilde E$ and every
$f\in \O(\D,X)\cap C(\overline \D,X)$ with $f(0)=x$ we have
$\sigma_f(\tilde E)=0$.
\end{theorem}

As a direct consequence of Theorems \ref{th:2} and \ref{th:3} we get
a characterization of pluripolar sets by analytic discs.

\begin{corollary}\label{cor:3}
Let $X$ be a Josefson manifold and assume that every bounded
plurisubharmonic function on $X$ is constant.
Let $E$ be a Borel subset in $X$. Then $E$ is pluripolar if and
only if
$$
\{x\in X\,;\, \exists f\in\O(\D,X)\cap C(\overline \D,X), f(0)=x, \sigma_f(E)>0\}
\neq X.
$$
\end{corollary}

Observe that even in $\C^n$ this corollary  gives a new
characterization of pluripolar sets.

\section{Analytic discs with images in boundaries of domains}

Let $X$ be a complex manifold and let $D\subset X$ be a domain.
If  $u\in\PSH(D)$, then we extend $u$
to an upper semicontinuous function  $u^*$ on the closure 
$\overline D$ by the formula
$$u^\ast(x)=\limsup_{D\ni y\to x}u(y), \qquad x\in \partial D.
$$
For every bounded function 
 $f:{\partial}D\to \R$  the function
\begin{equation*}
u_{f,D}=\sup\{v\in\PSH(D)\,;\, v^\ast|_{\partial D}\le f\}
\end{equation*}
is called the {\it Perron-Bremermann envelope} of $f$ on $D$.
We say that $D$ is {\it weakly regular} if for every relatively
open subset $U$ of ${\partial}D$ we have
$$
u_{-\chi_U,D}^*\leq -\chi_U \quad \text{ on } \ {\partial}D,
$$
where $\chi_U$ is the characteristic function of $U$.
We put $\omega(\cdot,U,D)=u_{-\chi_U,D}$.
Note that $\omega(\cdot,U,D)$ is a maximal plurisubharmonic function on
$D$ and $\omega^\ast(\cdot,U,D)\le-\chi_U$ on ${\partial}D$,
if $D$ is weakly regular.

We say that $D$ is {\it locally weakly regular} if for any $x\in
{\partial}D$ there exists a neighbourhood basis $\{V_j\}_{j=1}^\infty$
of $x$ in $X$ such that $D\cap V_j$ is weakly regular for all $j$.

Note that every locally weakly regular domain is weakly regular. Indeed,
for any $D_1\subset D_2$ and any open subset $U_j$ of   ${\partial}
D_j$ such that $U_1\subset U_2$ we have
$$
u_{-\chi_{U_1},D_1} \geq u_{-\chi_{U_2},D_2}.
$$

\begin{proposition}{\rm (cf. \cite{B2}).}
Any bounded domain $D$ in $\C^n$ which is regular for the Dirichet
problem ( as a domain in $\R^{2n}$) for the Laplace operator
is locally weakly regular.
In particular, any hyperconvex domain
is locally weakly regular.
\end{proposition}

\begin{proof}
The intersection of two Dirichlet regular domains is Dirichlet
regular, so it is enough to show that $D$ is weakly regular.
For a Dirichlet regular domain it is well-known that for any
$f\in C({\partial}D)$ we have $u_{f,D}^*\leq f$ on ${\partial}D$.
Since $-\chi_U$ is upper semicontinuous on ${\partial}D$ it is
sufficient to show that  $u_{f,D}^*\leq f$ on ${\partial}D$
for any upper semicontinuous function $f$.  Let $f_j$ be a sequence of
continuous functions decreasing to $f$.  Then
$u_{f,D}^*\leq u_{f_j,D}^*\leq f_j$ on ${\partial}D$.  We let
$j\to \infty$ and get $u_{f,D}^*\leq f$.
\end{proof}

\medskip

For any subset $A\subset\partial D$ we put
\begin{equation*}
\omega(x,A,D)=\sup\{\omega(x,U,D)\,;\,U \text{ is open and } 
A\subset U\subset\partial D\}, \qquad x\in D,
\end{equation*}
and
\begin{equation*}
\Omega(x,A,D)=-\sup\{\sigma_f(A)\,;\, 
f\in\O(\D,D)\cap C(\overline\D,\overline D), f(0)=x\},
\qquad x\in D.
\end{equation*}

We have a natural inequality between $\omega(\cdot,A,D)$ and
$\Omega(\cdot,A,D)$
as in the case when $A$ is in the interior of $D$.

\begin{lemma}
Let $X$ be a complex manifold,  $D\subset X$ be a weakly regular 
domain and $A\subset\partial D$ be a Borel set. Then
$\omega(\cdot,A,D)\le\Omega(\cdot,A,D)$.
\end{lemma}

\begin{proof} Let $x\in D$. If   $U$ is an open set in ${\partial}D$
such that $A\subset U\subset\partial D$,  and 
$f\in\O(\D,D)\cap C(\overline\D,\overline D)$  such that 
$f(0)=x$,  then for $u=\omega(\cdot,U,D)$ we have 
\begin{equation*}
u(x)\leq \int_\T u\circ f\, d\sigma
\leq
\int_{f^{-1}(A)\cap \T} u \circ f \, d\sigma
\leq -\sigma_f(A).
\end{equation*}
If we take supremum over $U$ in the left-hand side and infimum over $f$ 
in the right-hand side, then the inequality follows.
\end{proof}

\medskip
Now we will give a new proof of an improved version of
Lemma 9.1 in Poletsky \cite{P0}.

\begin{theorem}\label{corol:6} Let $X$ be a complex manifold, $D$ be a
relative compact weakly regular domain in $X$, and
 $U\subset\partial D$ be an open set. Then
$\omega(\cdot,U,D)=\Omega(\cdot,U,D)$.
\end{theorem}

The proof is in several steps each of which we state as a lemma.

\begin{lemma} Assume that $U_1\subset U_2\subset\dots\subset\partial D$
are open sets. Put $U=\cup_j U_j$. Then
\begin{equation*}
\lim_{j\to\infty}\omega(x,U_j,D)=\omega(x,U,D), \qquad x\in D.
\end{equation*}
\end{lemma}

\begin{proof} Put $u(x)=\lim_{j\to\infty}\omega(x,U_j,D)$ for 
$x\in D$. Note that the sequence is decreasing, so  
$u\in\PSH(D)$ and
$u\ge\omega(\cdot,U,D)$. On the other hand, 
$u^\ast\le\omega^\ast(\cdot,U_j,D)\le-\chi_{U_j}$ on
$\partial D$ for all $j\ge1$,  so $u^\ast\le-\chi_U$ on $\partial D$ 
and $u\le\omega(\cdot,U,D)$.
\end{proof}

\begin{lemma}\label{lemma:7}
For every $x_0\in U$ and  $\varepsilon>0$ there exists $r>0$ such that
\begin{equation*}
B(x_0,r)\cap D\subset\{x\in D\,;\,\omega(x,U,D)<-1+\varepsilon\}.
\end{equation*}
\end{lemma}

\begin{proof} Assume that for 
any $n\in\N$ there exists $x_n\in B(x_0,\frac 1n)\cap D$
such that $\omega(x_n,U,D)\ge-1+\varepsilon$. Then $x_n\to x_0$ and 
$\omega^\ast(x_0,U,D)\ge-1+\varepsilon$.
But $\omega^\ast(\cdot,U,D)\le-\chi_U$ on $\partial D$, a contradiction.
\end{proof}

\begin{lemma}\label{lemma:8}
Assume that $V\subset D$ is an open set such that for any 
$x_0\in U$ there exists an $r>0$ with
$B(x_0,r)\cap D\subset V$. Then
$\omega(\cdot,V,D)\le\omega(\cdot,U,D)$.
\end{lemma}

\begin{proof} We have 
$\omega^\ast(\cdot,V,D)\le-\chi_U$ on $\partial D$. Hence,
$\omega(\cdot,V,D)\le\omega(\cdot,U,D)$ on $D$.
\end{proof}

\begin{lemma}\label{lemma:5} We have
$\Omega(\cdot,\overline{U},D)\le\omega(\cdot, U,D)$.
\end{lemma}

\begin{proof} Fix $x_0\in D$ and $\varepsilon>0$ and take $a>\omega(x_0,U,D)$.
Let $\varepsilon_m>0$ be a sequence such that
$\varepsilon_m\searrow 0$. Assume that $\{B_j^m\}_{j}$ is a countable covering of $U$ with
balls of radii $<\varepsilon_m$ and centres in $U$ for any $m\ge1$. Put
\begin{equation*}
U_m=\cup_{j}\{x\in D\cap B_j^m\,;\,\omega(x,U\cap B_j^m, D\cap B_j^m)<-1+\varepsilon_m\}.
\end{equation*}
Let us show that
\begin{equation*}
U_m\subset\cup_{j}\{x\in D\cap B_j^m\,;\,\omega(x,U_{m+1}\cap B_j^m, D\cap B_j^m)<-1+\varepsilon_m\}.
\end{equation*}
For this, it suffices to show that
\begin{equation*}
\omega(\cdot,U_{m+1}\cap B_j^m, D\cap B_j^m)\le\omega(\cdot,U\cap B_j^m, D\cap B_j^m)\quad\text{ on }D\cap B_j^m,
\end{equation*}
for any $m\ge1$. It follows from Lemmas \ref{lemma:7} and  \ref{lemma:8}.

We have  $\omega(x_0,U,D)\ge\omega(x_0,U_1,D)$. So, there exists
an $f_1\in\O(\overline\D,D)$ such that $f_1(0)=x_0$ and
$\sigma_{f_1}(U_1)>|a|$. Put $\Delta_1=f_1^{-1}(U_1)\cap\T$. Now
we construct inductively $(f_m,\Delta_m)$ as in the proof of
Theorem~\ref{th:1} and get the analytic disc $f$ as its limit.
\end{proof}

\begin{prooftx} {Proof of Theorem~\ref{corol:6}}  By Lemma 4.2, we
have  
$\omega(\cdot,U,D)\le\Omega(\cdot,U,D)$.
Take open sets $U_1\Subset U_2\Subset\dots\Subset U$  such that $U=\cup_{j=1}^\infty U_j$.
According to Lemma \ref{lemma:5} we have $\Omega(\cdot,U,D)\le\omega(\cdot,U_m,D)$ on $D$  for any $m\ge1$.
Take $m\to\infty$ and get $\Omega(\cdot,U,D)\le\omega(\cdot,U,D)$ on $D$.
\end{prooftx}

\medskip
We let $X$ be a complex manifold, $D$ be a relatively compact
weakly regular domain in $X$, and  $A\subset\partial D$ be a Borel set.
Note that, in general, it is not true that $\{z\in A \,;\,
\omega^\ast(z,K,D)>-1\}$  is pluripolar.

\medskip
\begin{example} Let $K\subset\T$ be a non-polar compact set of
measure zero (take for example a Cantor-type
set on the unit circle). Then $\omega^\ast(\cdot,K,\D)$ is given as a
Poisson integral over the set $K$ which is of measure zero and
therefore $\omega^\ast(\cdot,K,\D)\equiv 0$.
\end{example}

\medskip
Nevertheless we have the following result, which is stated in
Sadullaev \cite{Sad}, Theorem~27.3.

\begin{lemma} Let $A\subset\partial D$ be a compact set.
Then $\omega^\ast(\cdot,A,D)\equiv0$ if and only if there exists
a $u\in\PSH^{-}(D)$, $u\not\equiv-\infty$,
such that $u^\ast|_{A}\equiv-\infty$.
\end{lemma}

\begin{proof}
If $\omega^*(\cdot,A,D)=0$, then there exists $x_0\in D$ such that
$\omega(x_0,A,D)=0$.  By the definition of $\omega(\cdot,A,D)$,
 there exists a sequence of
open sets $U_n\supset A$ on ${\partial}D$ so that
$\omega(x_0,U_n,D)>-2^{-n}$.  Now put $u=\sum_n\omega(\cdot,U_n,D)$.
Then $u\in \PSH(D)$, $u(x_0)>-1$, and $u^*|_A=-\infty$.

If, on the other hand, there exists a function $u\in \PSH(D)$ such
that $u^*|_A=-\infty$, then we consider the neighbourhoods
$U_n=\{u^*<-n\}\cap {\partial}D$ of $A$.  Then
$$
 \dfrac un \leq \omega(\cdot, U_n,D)\leq \omega(\cdot,A,D).
$$
So, $\omega(\cdot,A,D)=0$ on the set $\{u\neq -\infty\}$.
\end{proof}

\medskip
We say that a set $A\subset\partial D$ is {\it locally pluriregular
at $x_0\in \overline{A}$ with respect to $D$} if
there exists a sequence $r_j\searrow 0$ such that 
$\omega^\ast(x_0,B(x_0,r_j)\cap A,B(x_0,r_j)\cap D)=-1$ for any
$j\ge1$.

\medskip
\begin{theorem}\label{thm:10}
Let $A\subset\partial D$ be a Borel set.
Assume that $A=A_1\cup E$, where $A_1$ is locally pluriregular with
respect to $D$ and $E$ is
such that there exists a $u\in\PSH^{-}(D)$, $u\not\equiv-\infty$, with $u^\ast|_{E}\equiv-\infty$.
Then
\begin{equation*}
\omega(x,\overline A,D)\le\Omega(x,\overline A,D)\le\omega^\ast(x,A,D)\quad x\in D.
\end{equation*}
\end{theorem}

\begin{proof} Note that $\omega^\ast(\cdot,A,D)=\omega^\ast(\cdot,A_1,D)$ on $D$.
Fix $x_0\in D$ and $a>\omega^\ast(x_0,A,D)$.
Let $\varepsilon_m>0$ be a sequence such that
$\varepsilon_m\searrow 0$. Assume that $\{B_j^m\}_{j}$ is a countable
covering of $A_1$ with
balls of radii $<\varepsilon_m$ and centres in $A_1$ for any $m\ge 1$. Put
\begin{equation*}
U_m=\cup_{j}\{x\in D\cap B_j^m\,;\,\omega^\ast(x,A_1\cap B_j^m, D\cap B_j^m)<-1+\varepsilon_m\}.
\end{equation*}
Let us show that
\begin{equation*}
U_m\subset\cup_{j}\{x\in D\cap B_j^m\,;\,\omega(x,U_{m+1}\cap B_j^m, D\cap B_j^m)<-1+\varepsilon_m\}.
\end{equation*}
For this, it suffices to show that
\begin{equation*}
\omega(\cdot,U_{m+1}\cap B_j^m, D\cap B_j^m)\le\omega^\ast(\cdot,A_1\cap B_j^m, D\cap B_j^m)\quad\text{ on }D\cap B_j^m,
\end{equation*}
for any $m\ge1$. 

We have  $\omega^\ast(x_0,A_1,D)\ge\omega(x_0,U_1,D)$. So, there
exists an $f_1\in\O(\overline\D,D)$ such that $f_1(0)=x_0$ and
$\sigma_{f_1}(U_1)>|a|$. Put $\Delta_1=f_1^{-1}(U_1)\cap\T$. Now
we construct inductively $(f_m,\Delta_m)$ as in the proof of
Theorem~\ref{th:1} and get the analytic disc $f$ as its limit.
\end{proof}

\begin{question} Is Theorem \ref{thm:10} true for any Borel subset $A$
of $\partial D$?
\end{question}

\medskip\noindent
{\bf Acknowledgment.}
This paper was written while the first author was visiting the University of Iceland.
He likes to thank the Science Institute for its warm hospitality.

\bigskip
{\small
Armen Edigarian

Institute of Mathematics, Jagiellonian University, Reymonta 4,
30-059 Krak\'ow, Poland

E-mail: Armen.Edigarian@im.uj.edu.pl

\medskip

Ragnar Sigurdsson

Science Institute, University of Iceland, Dunhaga 3, IS-107
Reykjavik, Iceland

E-mail: ragnar@hi.is
}
\end{document}